\documentclass[12pt]{article}
\usepackage{amssymb,graphicx}

\def\d{{\rm d}}
\def\mi{{\rm i}}
\def\eps{\varepsilon}
\def\g{\gamma}
\def\G{\Gamma}

\def\l{\lambda}
\def\L{\Lambda}
\def\z{\zeta}
\def\Re{\mathop{\rm Re\,}\nolimits}
\def\Im{\mathop{\rm Im\,}\nolimits}
\def\e{\mathop{\rm e}\nolimits}

\def\hf{{\textstyle{1 \over 2}}}
\def\qt{{\textstyle{1 \over 4}}}
\def\tq{{\textstyle{3 \over 4}}}

\def\defi{\stackrel{\rm def}{=}}
\def\si{\!\!\! &}
\def\se{& \!\!\!}
\def\ni{\noindent}

\newcommand{\beq}{\begin{equation}}
\newcommand{\eeq}{\end{equation}}
\newcommand{\bea}{\begin{eqnarray}}
\newcommand{\eea}{\end{eqnarray}}

\title{Simplifications of the Keiper/Li approach to the Riemann Hypothesis}

\author{{\bf Andr\'e Voros}\\
CEA-Institut de Physique Th\'eorique de Saclay (CNRS-URA 2306)\\
F-91191 Gif-sur-Yvette Cedex (France)\\
E-mail: \tt andre.voros@cea.fr}

\begin{document}

\maketitle

\begin{abstract}
The Keiper/Li constants $\{ \l _n \}_{n=1,2,\ldots}$ are asymptotically ($n \to \infty$) sensitive 
to the Riemann Hypothesis, but highly elusive analytically and difficult to compute numerically.
We present quite explicit variant sequences that stay within the abstract Keiper--Li frame,
and appear simpler to analyze and compute.
\end{abstract}

The present work develops results that we announced in 2015.~\cite{VU}

\section{Generalities and notations}

We use the standard basic notions (e.g., \cite[chap. 8]{Da}):

\ni $\z(x):$ the Riemann zeta function (analytic over ${\mathbb C} \setminus \{+1\}$);

\ni $2 \xi (x):$ a completed zeta function, with its Riemann's Functional Equation:
\beq
\label{CZ}
2 \xi (x) \defi x(x-1) \pi^{-x/2}\G (x/2) \, \z (x) \equiv 2 \xi (1-x)
\eeq
(this \emph{doubled} Riemann's $\xi $-function is better normalized: $2 \xi (0) = 2 \xi (1) = 1$).

\ni $\{\rho \}:$ the set of zeros of $\xi $ (i.e., the nontrivial zeros of $\z $, or \emph{Riemann zeros},
counted with multiplicities if any, and grouped in pairs $(\rho,1-\rho )$ in the sums
that we write as $\sum_{\langle \rho,1-\rho \rangle}$); they all lie in the strip $\{ 0< \Re x <1 \}$.

\ni \textbf{Riemann Hypothesis (RH)} : \emph{all} the Riemann zeros lie on the \emph{critical line} 
$L \defi \{ \Re x = \hf  \}$.

\ni $k!! :$ the double factorial, to serve here for odd integers $k$ only, in which case
\bea
k!! \si\defi\se k(k-2) \cdots 1 \qquad \qquad \quad \mbox{for odd } k>0 , \nonumber \\
\si\defi\se 2^{(k+1)/2} \, \G (\hf k+1) / \sqrt \pi \quad \mbox{for odd } k\gtrless 0 \qquad (e.g.,\ (-1)!!=1) .
\eea
\ni $B_{2m}:$ the Bernoulli numbers; \quad $\g :$ Euler's constant.

\subsection{The Keiper and Li coefficients}

In 1992 Keiper \cite{K} considered a real sequence $\{ \l _n \}$ of generating function
\beq
\label{KEd}
f(z) \defi \log 2 \xi (M(z)) \equiv \sum_{n=1}^\infty \l ^{\rm K}_n \, z^n ,
\qquad M(z) \defi \frac{1}{1-z}, 
\eeq
($\l ^{\rm K}_n :$ our notation for \emph{Keiper's} $\l _n$), deduced that
\beq
\label{KEp}
\l ^{\rm K}_n \equiv n^{-1} \sum_{\langle \rho,1-\rho \rangle} \bigl[ 1-(1-1/\rho)^n \bigr] ,
\eeq
and argued that, under RH, $\l ^{\rm K}_n >0 \ (\forall n)$ and moreover
``if [...] the zeros are very evenly distributed, we can show that"  [this without proof]
\beq
\label{Kap}
\l ^{\rm K}_n \approx \hf \log n + c, \qquad c = \hf (\g -\log 2\pi -1) \approx -1.130330700754 \, .
\eeq

In (\ref{KEd}), the conformal mapping $M : x=(1-z)^{-1}$ acts to pull back 
the  critical line $L$ to the unit circle $\{ |z|=1 \}$, with the fundamental consequence:
\beq
\label{KLS}
\mbox{RH } \iff \ f \mbox{ \emph{regular in the whole open unit disk} } \{ |z|<1 \}.
\eeq
Then, (\ref{KEd}) specifies the sequence $\{ \l _n \}$ as a particular encoding 
of the germ of $\log 2\xi (x)$ at the ``basepoint" $\underline x=M(0)$ (here:~$\underline x=1$).
\medskip

In 1997 Li \cite{LI1} independently introduced another sequence $\l _n$, through
\beq
\label{LId}
\l ^{\rm L}_n = \frac{1}{(n-1)!}\, \frac{\d ^n}{\d x^n} \bigl[ x^{n-1} \log 2 \xi (x) \bigr]_{x=1},
\quad n=1,2,\ldots \quad (\l ^{\rm L}_n \defi \mbox{ \emph{Li's} } \l _n ),
\eeq
deduced that
\beq
\label{LIp}
\l ^{\rm L}_n \equiv \sum_{\langle \rho,1-\rho \rangle} \bigl[ 1-(1-1/\rho)^n \bigr] ,
\eeq
and proved the \emph{sharp equivalence}: RH $\iff \ \l ^{\rm L}_n >0$ for all $n$ (\emph{Li's criterion}).

Actually, by comparing (\ref{KEp}) vs (\ref{LIp}) for instance,
\beq
\l ^{\rm L}_n \equiv n \, \l ^{\rm K}_n \qquad \mbox{for all } n=1,2,\ldots ;
\eeq
our superscripts K vs L will disambiguate $\l _n$ whenever the factor $n$ matters.

\subsection{Probing RH through the Keiper--Li $\{ \l _n \}$}

In 2000 Oesterl\'e proved (but left unpublished) \cite{O} that RH alone implies
\beq
\label{OE}
\l ^{\rm L}_n = n (\hf \log n + c) + o(n), \quad 
\mbox{with } c = \hf (\g -\log 2\pi -1) \, \mbox{ as in (\ref{Kap})} .
\eeq

In 2004--2006, using the saddle-point method on an integral form of $\l _n$,
we gave an \emph{asymptotic criterion} for RH \cite{V}\cite{VB} in the form of this alternative: 
\bea
\label{ASF}
\bullet \ \mbox{RH false:} \quad \l ^{\rm L}_n \si\sim\se - \!\sum_{\Re \rho >1/2}(1-1/\rho)^{-n} 
\pmod {o(r^n)\ \forall r>1} ; \\
\label{ASC}
\bullet \ \mbox{RH true:} \quad \l ^{\rm L}_n \si\sim\se n (\hf \log n + c) \pmod {o(n)}
\eea
(erratum: we had the sign wrong in the case RH false, which did not affect the purely qualitative consequences 
we drew at the time; correction in \cite{VK}).

In 2007 Lagarias \cite{Lg} strengthened (\ref{OE}) by improving $o(n)$ to $O(\sqrt n \log n)$.
\medskip

To assess how the above criteria may advance the testing of RH, 
one must bring in the height $T_0$ up to which RH is confirmed by direct means:
\beq
T_0 \approx 2.4 \cdot 10^{12} \quad \mbox{currently (since 2004). } \cite{G}
\eeq
It is then known that: first, no $\l _n$ can go negative as long as $n<T_0^{\, 2}$ \cite{O}\cite{BPY};
and more broadly, if a zero $\rho = \hf \pm t \pm \mi T$ violates RH (with $t>0$, $T >T_0$),
then no effect of that will be detectable upon the $\l _n$ unless \cite{VB}
\beq
\label{LIt}
n \gtrsim T^2/t > 2T_0^{\, 2} \quad \mbox{(since $t<\hf $), currently implying} \quad n \gtrsim 10^{25} 
\eeq
($n \gtrsim T^2/t$ actually states the \emph{uncertainty principle} in the relevant geometry).

At the same time, the $\l _n$ are quite elusive analytically \cite{BL}\cite{C2},
and also numerically (see Ma\'slanka \cite{M1}\cite{M2} and Coffey \cite{C1})
as their evaluation requires a recursive machinery, whose intricacy grows very rapidly with $n$,
and which moreover destroys ca. $\qt $ decimal place of precision per step~$n$ 
(if done \emph{ex nihilo} - i.e., using no Riemann zeros as input) \cite[fig.~6]{M2}.
Thus only $\l _n$-values up to $n \approx 4000$ have been accessed \emph{ex nihilo}, 
so that the useful range (\ref{LIt}) looks way beyond reach.

\section{An \emph{explicit} variant sequence $\{ \L _n \}$}

We propose to deform the $\{ \l ^{\rm K}_n \}$ (in Keiper's normalization (\ref{KEd}))  
into a simpler sequence $\{ \L _n \}$ having a \emph{totally closed form}. 
The original $\l _n$ appeared rigidly specified, but only inasmuch as the pole $x=1$ of $\z (x)$
was invariably made \emph{the} basepoint. 
Now while this choice can make sense, it is by no means compulsory.
On the contrary, other conformal mappings than $M$ in (\ref{KEd}) 
realize the Keiper--Li \emph{idea} (RH-sensitivity, embodied in (\ref{KLS})) just as well: 
the key condition is that all Riemann zeros on $L$ must pull back to $\{ |z|=1 \}$, 
achieving (\ref{KLS}),
while nothing binds the basepoints $\underline x$ to which $z=0$ can map;
the resulting $\l _n$ will just vary with $\underline x$ 
as functions of the derivatives $\xi ^{(m)}(\underline x)$.
As such, Sekatskii's ``generalized Li's sums" \cite{Se} have 
$\underline x = (1-a) \in {\mathbb R} \setminus \{ \hf \}$,
whereas our ``centered" $\l ^0_n$ were tailored to have $\underline x = \hf $, 
the symmetry center for $\xi (x)$ (\cite[\S~3.4]{VK}, and Appendix).
Our next construction will push this idea of deformation even further,
and have no single distinguished basepoint (except, loosely, $\underline x=\infty \,$?):
we will substitute selected \emph{finite differences}
for the derivatives of $\log \xi$ that enter the original $\l _n$ 
(and, in the Appendix, our centered $\l^{(0)}_n$), to attain more explicit sequences.

\subsection{Construction of the new sequence}

The original definition (\ref{KEd}) is equivalent, by the residue theorem, to the contour integral formula
\beq
\label{Kint}
\l ^{\rm K}_n = \frac{1}{2 \pi\mi} \oint \frac{\d z}{z^{n+1}} f(z) , 
\qquad f(z) \equiv \log 2 \xi \Bigl(\frac{1}{1-z}\Bigr) ,
\eeq
with a positive contour in the unit disk around $z=0$ excluding all other singularities 
(i.e., those of $f$). Derivatives of $\log 2\xi (x)$ up to order $n$ occur in $\l _n$
because the denominator $z^{n+1}$ has all its zeros degenerate (at $z=0$). 

Now at given $n$, if we \emph{split those zeros apart} as $0, z_1, \ldots, z_n$ 
(all distinct, and still inside the contour), 
then the so modified integral evaluates to a linear combination of the $f(z_m) :$ 
derivatives become \emph{finite differences}.
To split the zeros, instead of plain shifts of the factors $z \mapsto z-z_m$ 
which fail to preserve the all-important unit disk, we use \emph{hyperbolic} translates 
\beq
{z \mapsto B_{z_m}(z)} = {(z-z_m)/(1-z_m^\ast z)} \quad \mbox{(M\"obius transformations)}.
\eeq
The point $z=0$ has now lost its special status, hence so does the particular mapping $M$ 
(picked for pulling back the pole $x=1$ to $z=0$),
so that the variable $x$, natural for the $\z $-function, also becomes the simplest to use. 
Then (\ref{Kint}) expresses as
\beq
\label{KLX}
\l ^{\rm K}_n = \frac{1}{2 \pi\mi} \oint \frac{\d x}{x(x-1)}
\Bigl( \frac{x}{x-1} \Bigl)^n \log 2\xi (x)  \quad ( \mbox{integrated around }x=1)  ,
\eeq
and the deformations as above read as
\beq
\label{KLD}
\frac{1}{2 \pi\mi} \oint_{{\cal C}_n} \frac{\d x}{x(x-1)} 
\frac{1}{b_{x_1}(x) \ldots b_{x_n}(x)} \log 2\xi (x) , \ \quad 
b_{\tilde x} (x) \equiv \frac{\tilde x^\ast }{\tilde x} \, \frac{x-\tilde x}{x + \tilde x^\ast -1},
\eeq
where the contour ${\cal C}_n$ encircles the points $1,x_1,\ldots,x_n$ positively
(and may as well depend on~$n$). Now the integral in (\ref{KLD}) readily evaluates to
\beq
\label{Bint}
\sum_{m=1}^n \frac{1}{x_m (x_m-1)} \, \frac{1}{[ b_{x_1} \ldots b_{x_n} ]' (x_m) } \log 2\xi (x_m)
\eeq
by the residue theorem ($x=1$ contributes zero since $\log 2\xi (1) =0$).

Finally, for each $n$ we select $x_m \equiv 2m$ for $m=1,2,\ldots$ (independently of~$n$) 
to benefit from the known values $\z (2m)$, 
and a contour ${\cal C}_n $ just encircling the real interval $[1,2n]$ positively
(encircling the subinterval $[2,2n]$ would suffice, 
however here it will always be of interest to dilate, not shrink, ${\cal C}_n $). 
All that fixes the sequence
\bea
\label{Ldef}
\L _n \si\defi\se \frac{1}{2 \pi\mi} \oint_{{\cal C}_n } \frac{\d x}{x(x - 1)} \, G_n(x) \log 2 \xi (x) , \\
\label{GN}
G_n(x) \si\defi\se \prod_{m=1}^n \frac{x + 2m - 1}{x-2m} \equiv 
\frac{ \G (\hf x-n) \, \G (\hf (x \!+\! 1) + n) }{ \G (\hf x) \, \G (\hf (x \!+\! 1)) } \\
\label{gdf}
\si\equiv\se g(x) (-1)^n \frac{\G (\hf (x \!+\! 1) + n)}{\G(1-\hf x + n)}, \qquad
g(x) \defi \frac{\sqrt \pi \, 2^{x-1}}{\sin (\pi x/2) \, \G (x)} 
\eea
(by the duplication and reflection formulae for $\G $). For this case, (\ref{Bint}) yields
\beq
\label{EKL}
\L _n \equiv (-1)^n \sum_{m=1}^n (-1)^m A_{nm} \log 2 \xi (2m) , \qquad \qquad n=1,2,\ldots ,
\eeq
with
\bea
\label{Anm}
A_{nm} \si=\se \frac{2^{-2n}}{2m-1} {2(n \!+\! m) \choose n \!+\! m} {n \!+\! m \choose 2m} 
\equiv \frac{2^{m-n} \, \bigl( 2(n+m)-1 \bigr) !!}{(2m-1) \, (n-m)! \, (2m)!} \qquad \nonumber \\
\si\equiv\se \frac{2^{2m} \G (n+m+1/2)} {(2m-1) \, (n-m)! \, (2m)! \, \sqrt\pi}
\qquad \qquad (\mbox{for } m=0,1,2,\ldots) , \\
\label{BE}
2 \xi (2m) \si=\se \frac{ |B_{2m}| }{ |(2m-3)!!| } \, (2\pi)^m \equiv 
\frac{ 2 \, |B_ {2m}| }{ |\G (m-\hf )| } \, \pi^{m+1/2}
\eea
(the absolute values in the last two denominators only act for $m=0$, 
resulting in $\log 2 \xi (0)=0$ which vanishes thereafter).

So, this particular deformation $\{ \L _n \}$ of Keiper's $\{ \l ^{\rm K}_n \}$ is specified by (\ref{EKL}) 
in a \emph{totally explicit} form (and fairly uniquely dictated as above).
With no recursion involved, any single $\L _n$ can be computed straight away and by itself, 
in welcome contrast to the original $\l _n$. 
\medskip

\textbf{Remarks.}
\smallskip

\ni 1) $\sum\limits_{m=1}^n (-1)^m A_{nm} m$ is computable by the second sum rule (\ref{AID}) below
(with $A_{n0} \equiv  - 2^{-n} (2n-1)!! \, / \, n!$ by (\ref{Anm}));
the  $(\log 2\pi)$-contributions to (\ref{EKL}) from the first expression (\ref{BE})
can thereby be summed, resulting in $\L_n \equiv \hf \log 2\pi + u_n $ with
\beq
\label{Udef}
u_n \defi (-1)^n \left[ \sum_{m=1}^n (-1)^m A_{nm} \log \frac{ |B_{2m}|}{(2m-3)!!} 
+ \frac{1}{2 A_{n0}} \log 2\pi \right] :
\eeq
it was through this sequence $\{ u_n \}$ that we earlier announced our results \cite{VU}.
Likewise, the last expression (\ref{BE}) leads to the partially summed form
\bea
\label{Vdef}
\L _n \si\equiv\se \hf \log \pi + (-1)^n 
\left[ \sum_{m=1}^n (-1)^m A_{nm} \log \frac{ |B_{2m}|}{\G (m \!-\! \hf )} \right. \nonumber \\
\si\se \qquad \qquad \qquad \quad \ \left. + \Bigl( \frac{1}{A_{n0}}\!-\! A_{n0} \Bigr) \log 2 
+ \Bigl( \frac{1}{A_{n0}}\!-\! \frac{A_{n0}}{2} \Bigr) \log \pi \right] .
\eea

\smallskip

\ni 2) if in place of (\ref{BE}) we use (\ref{CZ}) and the expanded logarithm of the Euler product:
$\log \z (x) \equiv \sum\limits_p \sum\limits_{r=1}^\infty p^{-rx}/r \quad (x>1)$
where $p$ runs over the \emph{primes}, 
then (\ref{EKL}) yields an \emph{arithmetic} form for $\L _n$, in analogy with \cite[thm~2]{BL} for $\l ^{\rm L}_n$.
\smallskip

\ni 3) B\'aez-Duarte's sequential criterion for RH \cite{BD} 
is similarly explicit in terms of the Bernoulli numbers,
but there, any effect of RH-violating zeros seems hopelessly tiny until inordinately large $n \gtrsim \e ^{\pi T_0}$
\cite[\S~4]{M3}\cite[\S~7]{FV} (the latter quotes $n \gtrsim 10^{600,000,000}$).

\ni 4) With L-functions for real primitive Dirichlet characters $\chi $ in place of $\z $, 
\cite[chaps.~5, 6, 9]{Da}
the whole argument carries over, essentially unchanged for $\chi $ even, whereas
\beq
\L _{\chi ,n} = (-1)^n \sum_{m=1}^n (-1)^m 
\frac{2^{m-n} \, \bigl( 2(n+m)+1 \bigr) !!}{(2m \!+\! 1) \, (n \!-\! m)! \, (2m \!+\! 1)!} 
\log \xi_\chi (2m+1)
\eeq
for $\chi $ odd, where $\xi _\chi (x)$ is the completed $L _\chi $-function 
(normalized to $\xi _\chi (0) \equiv {\xi _\chi (1) =1}$, like $2 \xi (x)$ for $\z $ in (\ref{CZ})), 
whose values at $x=2m+1$ are explicit.

\subsection{Expression of $\L _n$ in terms of the Riemann zeros}

Let the primitive
\bea
\label{Fdef}
F_n(x) \si\defi\se \int_\infty^x \frac{G_n(y)}{y(y-1)} \, \d y \nonumber \\
\si\equiv\se (-1)^n \biggl[ - \frac{1}{A_{n0}} \log (x-1)
+ \sum_{m=0}^n (-1)^m A_{nm} \log (x-2m) \biggr] 
\eea
be defined as single-valued from a neighborhood of $x=\infty$ 
to the whole $x$-plane minus the cut $[0,2n]$.
E.g., $F_1(x) = {\hf \log \, [x(x-2)^3/(x-1)^4]}$; 
and for general~$n$, (\ref{Fdef}) follows from, e.g., \cite[\S 2.102]{GR}
using the $A_{nm}$ from (\ref{Anm}).

For $x \to \infty, \ F_n(x) \sim \int_\infty^x \d y/y^2 = -1/x$;
the consistency of this with (\ref{Fdef}) imposes the identities
\beq
\label{AID}
\sum_{m=0}^n (-1)^m A_{nm} \equiv  \frac{1}{A_{n0}} , \qquad
2 \sum_{m=1}^n (-1)^m A_{nm} m \equiv  (-1)^n + \frac{1}{A_{n0}} .
\eeq

In terms of (\ref{Fdef}), the $\L _n$ result by summing the following series over the zeros 
(converging like $\sum_{\langle \rho,1-\rho \rangle} 1/\rho $ for any $n$):
\beq
\label{LIR}
\L _n \equiv \sum_{\langle \rho,1-\rho \rangle} F_n (\rho) , \qquad n=1,2,\ldots .
\eeq
(For the original $\l ^{\rm K}_n $, (\ref{Fdef}) uses $[x/(x-1)]^n$ in place of $G_n$ by (\ref{KLX}), 
exceptionally yielding \emph{rational} functions: ${n^{-1}[1-(1-1/(1-x))^n]}$, 
for which (\ref{LIR}) restores~(\ref{KEp}).)

Proof of (\ref{LIR}) (condensed, see fig. 1): first stretch the contour ${\cal C}_n$ in (\ref{Ldef}) 
to ${\cal C}_n '$ fully enclosing the cut $[0,2n]$ of $F_n$ (as allowed by $\log 2\xi (0)=0$). 
Since $F_n$ is single-valued on ${\cal C}_n '$, the so modified (\ref{Ldef}) can be integrated by parts,
\beq
\label{Lfed}
\L _n \defi -\frac{1}{2 \pi\mi} \oint_{{\mathcal C}_n '} F_n(x) \Biggl[ \frac{\xi '}{\xi } \Biggr] (x) \d x ,
\eeq
then the contour ${\mathcal C}_n '$ can be further deformed into a sum of 
an outer anticlockwise circle ${\mathcal C}_R$ centered at $\hf $ of radius $R \to \infty$ (not drawn),
and of small clockwise circles around the poles of the meromorphic function $\xi '/\xi $ 
inside ${\mathcal C}_R$; 
these poles are the Riemann zeros $\rho $, and each contributes $F_n(\rho )$.
By the Functional Equation (\ref{CZ}), the integral on ${\mathcal C}_R$ is also
$\oint_{{\mathcal C}_R} \!{\hf [F_n (x)+F_n (1 \!-\! x)]} \, [ \frac{\xi '}{\xi } ](x) \d x$, 
which tends to 0 if $R \to \infty$ staying far enough from ordinates of Riemann zeros in a classic way 
(so that $|\z '/\z | (s + \mi R) < K \log ^2 R$ for all $s \in [-1,+2]$ \cite[p.~108]{Da}), 
hence (\ref{LIR}) results. \hfill $\square$

\begin{figure}
\center
\includegraphics[scale=.6]{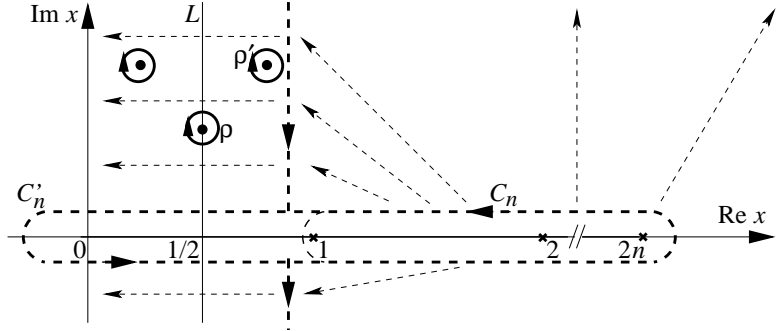}
\caption{\small Deformation of the integration path for the integral (\ref{Lfed}) 
against the meromorphic function $\xi '/\xi $
whose poles are the Riemann zeros, here exemplified - not on scale - by $\rho $ (on the critical line), 
and $\rho '$ (off the line, putative, shown with its partner across the critical line). 
A symmetrical lower half-plane is implied.}
\end{figure}

\section{Criterion for RH based on the new sequence}

We will sketch why the totally explicit sequence $\{ \L _n \}$ largely shares the sensitivity to RH
of the highly elusive Keiper--Li sequence.

\subsection{Asymptotic criterion}

We will mainly argue an \emph{asymptotic} sensitivity to RH as $n \to \infty$,
through this alternative for $\{ \L _n \}$ which parallels (\ref{ASF})--(\ref{ASC}) for $\{ \l _n \} :$
\bea
\label{LNRH}
\bullet \ \mbox{RH false:} \ \L _n \si\sim\se 
\Biggl\{ \sum_{\Re \rho >1/2} F_n (\rho ) \Biggr\} \pmod {o(n^\eps )\ \forall \eps >0} , \\
\label{ST1}
\mbox{and} \quad F_n (\rho ) \si\sim\se 
\frac{g(\rho )}{\rho (\rho -1)} (-1)^n \frac{n^{\rho -1/2}}{\log n} \qquad \qquad (n \to \infty) , \\
\label{LABS}
\Longrightarrow \quad |F_n(\rho )| \si\approx\se 
\frac{1}{|\Im \rho |^2 \log n} \Bigl( \frac{2n}{|\Im \rho |} \Bigr) ^{\Re \rho } 
\qquad \mbox{for } n \gg |\Im \rho | \gg 1 . \\[4pt]
\label{LRH}
\bullet \ \mbox{RH true:} \ \L _n \si\sim\se \log n + C , \quad C = \hf(\g \!-\! \log \pi \!-\! 1) 
\approx -0.783757110474 , \qquad
\eea
the latter to be compared to (\ref{OE}), with $C=c+\hf \log 2$. 
As for (\ref{LNRH}), the summation converges if the terms with $\rho $ and $\rho^\ast $ 
are grouped together (as symbolized by the curly brackets), and more caveats are issued in \S~3.2.
\medskip

We give a condensed derivation. Past some common generalities, 
we will separate the cases RH true/false (short of a unified method as in \cite{V}).
\smallskip

The general idea is nowadays known as large-order perturbative analysis
or instanton calculus, but initially we just follow the pioneering Darboux's theorem \cite[\S 7.2]{Di}\cite{BPV} 
to get the large-order behavior of Taylor series like (\ref{KEd}) out of the integral form (\ref{Kint})
or more simply, its integration by parts $\l ^{\rm L}_n = (2 \pi\mi)^{-1} \oint z^{-n} f'(z) \, \d z$
because $f'$ is meromorphic whereas $f$ has branch cuts.
Then this integrand has the large-$n$ form $\e^{\Phi_n (z)}$ where $\Phi_n $ tends to $\infty$ with $n$ 
($\Phi _n (z) \sim -n \log z$), hence the steepest-descent method applies: \cite[\S~2.5]{Er}
we deform the integration contour $C$ toward decreasing $\Re \Phi _n$, 
i.e. here, into a circle of radius growing toward 1 (fig.~2);
then, each of the encountered singularities of $f'$, 
here simple poles $M^{-1}(\rho ')$ for RH-violating zeros $\rho '$,
yields an asymptotic contribution $-z_{\rho '}^{\, -n}$, all of which add up to (\ref{ASF}). \cite{V}
If on the other hand RH is true, then the contour can arbitrarily approach the unit circle, 
(\ref{ASF}) stays empty, and only a finer analysis of the limiting integral 
(\cite{O}, recalled in \S~3.3.1 below) leads to a definite asymptotic form, as (\ref{OE}).

\begin{figure}[h]
\center
\includegraphics[scale=.6]{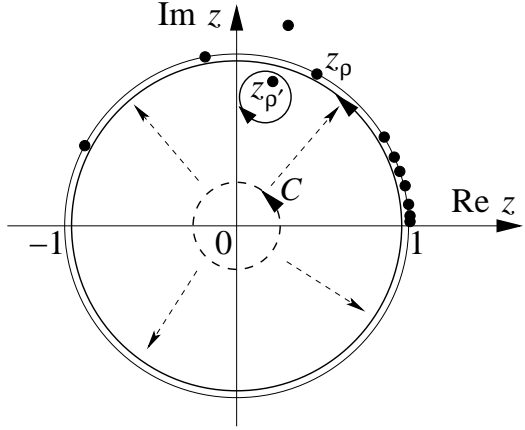}
\caption{\small Contour deformation for $\l ^{\rm L}_n, \ n \to \infty$, in the steepest-descent method.
(The symmetrical poles in the lower half-plane are not drawn.)
}
\end{figure}
\medskip

We then wish to do the same with an ($x$-plane) integral form for $\L _n$, 
be it (\ref{Ldef}) (with the function $G_n (x)$ defined by (\ref{GN})--(\ref{gdf})), 
or (\ref{Lfed}) (with $F_n (x)$ defined by (\ref{Fdef})). Now (\ref{gdf}) at once implies
\beq
G_n (x) \sim g(x) (-1)^n n^{x-1/2} \sim g(x) (-1)^n \e^{\log n (x-1/2)} \quad
\mbox{for } n \to \infty \mbox{ at fixed } x ,
\eeq
hence now the large asymptotic parameter is $\log n$ and the large-$n$ level lines
of the integrand are $\{ \Re x = \hf + t_0 \}$. 
For the steepest-descent method, $|z| \to 1^-$ in fig.~2 thus becomes $t_0 \to 0^+$.
A new complication is that these level lines now all terminate at $\infty$, an essential singularity.
Temporarily ignoring this, we note that the contour deformation on (\ref{Lfed}) for $\L _n$ 
has already yielded (\ref{LIR}), so we simply have to extract the asymptotically relevant part 
of $\sum_\rho F_n(\rho )$.
For $n \to \infty$, $F_n (\rho )$ is to be expressed using a steepest-descent path \cite[\S~2.5]{Er}, as
\beq
F_n (\rho ) = \int_{-\infty}^\rho \frac{G_n(x)}{x(x-1)} \, \d x
\sim \int_{-\infty}^\rho \frac{g(x)}{x(x-1)} (-1)^n n^{x-1/2} \, \d x :
\eeq
a Laplace transform in the variable $\log n$, of asymptotic form \cite[eq.~2.2(2)]{Er}
\beq
\label{FAS}
F_n (\rho ) \sim \frac{g(\rho )}{\rho (\rho -1)} (-1)^n \frac{n^{\rho -1/2}}{\log n} .
\eeq
Consequently, the removal of all $o(n^{t_0 + \eps} )$ terms from (\ref{LIR}) 
unconditionally leaves us with
\beq
\label{LSM}
\L _n = \Biggl\{ \sum_{\Re \rho >t_0} F_n (\rho ) \Biggr\} + o(n^{t_1}) \qquad
\mbox{for all } t_1>t_0 \ge 0 ,
\eeq
where the summation converges if the terms with $\rho $ and $\rho^\ast $ are grouped together
(as symbolized by the curly brackets).

However, in the RH true case, (\ref{LSM}) with $t_0=0$ delivers no better than
$\L _n = o(n^\eps )\ \forall \eps >0$,
and only a finer analysis of the limiting integral on the critical line $L$
will lead to a definite asymptotic form, in \S~3.3.2. Hence we pursue the case RH false first.

\subsection{Details for the case RH false}

If RH-violating zeros exist, we cannot ensure that they are finitely many, 
nor that they otherwise can be enumerated according to non-increasing real parts.
Then, unlike (\ref{ASF}), the series (\ref{LSM}) ought not to be directly readable 
as an explicit asymptotic expansion, to whatever order $t_0<\hf $. 
Instead, any closed-form asymptotic statement on $\L _n$ would have to involve
the detailed \emph{2D} distribution of RH-violating zeros toward~$\infty$, currently unknown.
Indeed, for no $t_0 < \hf $ can we perform or describe the sum of the series (\ref{LSM}) explicitly 
(barring the purely hypothetical case of finitely many terms).
In particular, it ought to be unlawful to substitute the individual asymptotic forms (\ref{FAS}) 
in bulk into any of the series (\ref{LSM}); we can only interpret the latter as a total 
of individual RH-violating zeros' contributions to the large-$n$ behavior of $\L _n$. 

Moreover, any such zero $\rho = \hf + t \pm \mi T$ with $t>0$ must have $T>T_0$ hence $T \gg 1$, 
which implies
\beq
|g(\rho )| \approx \Bigl( \frac{2}{T} \Bigr) ^t \quad \Longrightarrow \quad 
|F_n(\rho )| \approx \frac{1}{T^2 \log n} \Bigl( \frac{2n}{T} \Bigr) ^t .
\eeq
All in all, letting $t_0=0$ we obtain (\ref{LNRH})--(\ref{LABS}) in the case RH false.
\hfill $\square$
\medskip

The upshot of (\ref{ST1}) is that each RH-violating zero $\rho $ imparts $\L _n$
with a growing $n^{\rho -1/2}$-like oscillation; one consequence (in view of \S~3.4 below)
is that it would take improbable cancellations to have $\L _n >0$ forever, if RH was false.

\subsection{Details for the case RH true}

Here our quickest path is to adapt:

\subsubsection{Oesterl\'e's argument for $\l ^{\rm L}_n$}
(as reworded by us). We start from this real integral giving $\l ^{\rm K}_n :$ \cite{O}\cite{V}
\beq
\label{LR}
\l ^{\rm K}_n = \int_0^\pi 2 \sin n \theta \ N(\hf \cot(\hf \theta)) \, \d \theta \, ;
\eeq
here $N(T) = \# \{ \rho \in [\hf,\hf+\mi T] \subset L\}$ is the zeros' staircase counting function;
$T \equiv \hf \cot(\hf \theta )$ where $\theta \in(0,\pi]$ is the angle subtended 
by the real segment $[0,1]$ from the point $\hf+\mi T$,
$\d T \equiv -(\qt + T^2) \, \d \theta $, and the integrand is actually the reduction of
\beq
\label{LIX}
2 \Im \Bigl( \frac{x}{x-1} \Bigr)^n \log 2\xi (x) \frac{\d x}{x(x-1)}
\eeq
once the integration path in (\ref{Ldef}) has reached $\{x=\hf + 0 +\mi T\}$ (under RH)
and $\theta $ reparametrizes $T$.

Then $\l ^{\rm K}_n$  mod $o(1)$ will stem from the Riemann--von Mangoldt theorem: 
\cite[chaps.~8, 15]{Da}
\beq
N(T) = \frac{T}{2\pi} \Bigl( \log \frac{T}{2\pi}-1 \Bigr) + \delta N(T) , \quad
\delta N(T) = O(\log T) \, \mbox{ as } T\to +\infty.
\eeq
Proof: (\ref{LR}) mod $o(1)$ evaluates as follows: 

1) in $N(\cdot )$, the term $\delta N(\cdot )$ is integrable \emph{up to $\theta =0$ included},
then its integral against $\sin n \theta $ is $o(1)$ (Riemann--Lebesgue lemma) hence negligible;

2) change to the variable $\Theta _n \equiv n \theta $; then, change the resulting upper integration bound 
$n \pi$ to $+\infty$ and use $T \sim 1/\theta = n/\Theta _n$ to get, mod $o(1)$,
\beq
\l ^{\rm K}_n \sim \int_0^\infty 2 \sin \Theta _n \, \frac{n}{2 \pi \Theta _n} 
\Bigl[ \log \frac{n}{2 \pi \Theta _n} -1 \Bigr] \frac{\d \Theta _n}{n} .
\eeq
Now the classic formulae $\int_0^\infty \sin \Theta \, \d \Theta /\Theta = \pi /2$ 
and $\int_0^\infty \sin \Theta  \log \Theta \, \d \Theta / \Theta  = - \pi \g /2$ 
\cite[eqs.~(3.721(1)) and (4.421(1))]{GR} yield the result (amounting to (\ref{OE}))
\beq
\label{Kas}
\l ^{\rm K}_n = \hf \log n + c + o(1) \qquad \mbox{under RH true}.
\eeq
\hfill $\square$

\subsubsection{Parallel treatment for $\L _n$}

Basically for $\L _n$, $\displaystyle \Bigl( \frac{x}{x-1} \Bigr)^n$ in (\ref{LIX}) 
is to be replaced by $G_n(x)$ from (\ref{GN}), hence (\ref{LR}) changes to
\beq
\label{GR}
\L _n = \int_0^\pi 2 \sin \Theta _n (\theta) \ N(\hf \cot(\hf \theta)) \, \d \theta \, ,
\eeq
where $\Theta _n \in (0,n\pi]$ (previously $\equiv n \theta $) is now the sum of the $n$ angles 
subtended by the real segments $[1-2m,2m]$ from the point $\hf+\mi T$, for $m=1,2,\ldots,n$. 
The two \emph{endpoint slopes} of the function $\Theta _n (\theta )$ will mainly matter (independently):
\bea
\Theta _n'(0) \si=\se \sum_{m=1}^n (4m-1) \equiv n(2n+1) , \\
\Theta _n'(\pi) \si=\se \sum_{m=1}^n (4m-1)^{-1} \equiv 
\qt \bigl[ (\G '/\G )(n+\tq ) + \g + 3 \log 2 - \pi/2 \bigr] .
\eea

We then follow the same steps as with $\l ^{\rm K}_n$ just above.

1) $\int_0^\pi 2 \sin \Theta _n (\theta) \ \delta N(\hf \cot(\hf \theta)) \, \d \theta = o(1)$
if a nonstationary-phase principle can apply for the oscillatory function $\sin \Theta _n (\theta)$,
i.e., if the minimum slope of $\Theta _n (\theta)\ (\theta \in [0,\pi ])$ goes to $\infty$ with $n$:
previously (with $\Theta _n \equiv n \theta)$ that slope was $n$, 
now it is $\Theta _n' (\pi) \sim \qt \log n$ which still diverges for $n \to \infty$
therefore gives the $o(1)$ bound; but due to $\Theta _n' (\pi) \ll n$, 
this $o(1)$ may decay much slower than the corresponding $o(1)$ for $\l ^{\rm K}_n$.

2) In this step (i.e., $T \to +\infty$), only $\theta \to 0$ behaviors enter;
here $\Theta _n \sim \Theta _n'(0) \, \theta$, vs $n \theta$ previously, so it suffices to substitute 
$\Theta _n'(0)$ for $n$ in the asymptotic result (\ref{Kas}) for $\l ^{\rm K}_n$, to get
\beq
\L _n \sim \hf \log \Theta _n'(0) + c = \hf \log [n(2n+1)] + c \sim \log n + (c + \hf \log 2).
\eeq
\hfill $\square$

\subsection{Asymptotic or full-fledged Li's criterion?}

We do not control well enough the function $F_n$ in (\ref{Fdef}) or for that matter, the primitive 
$\int \sin \Theta _n (\theta ) \, \d \theta $ in (\ref{GR}), to be able to infer that RH implies $\L _n >0$ for all $n$,
as was the case for $\l _n$ straightforwardly from (\ref{KEp}).

On the other hand, our criterion (\ref{LNRH})--(\ref{LRH}) is synonymous
of large-$n$ positivity for $\L _n$ if and only if RH holds (invoking the last sentence of \S~3.2),
while low-$n$ positivity is numerically patent (see next \S).

All in all, as an aside we then also \emph{conjecture} that: 
\emph{Li's criterion works for the sequence $\L _n$} (RH $\ \iff \ \L _n >0$ for all $n$).

\section{Quantitative aspects}

\subsection{Numerical data}

\begin{figure}
\center
\includegraphics[scale=.4]{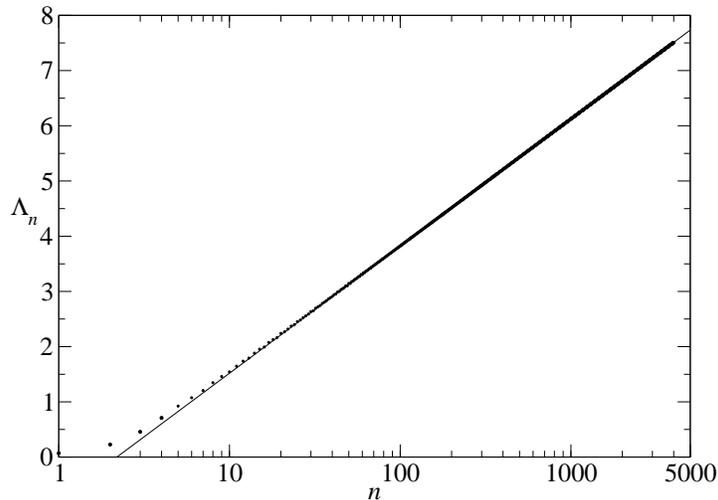}
\caption{\small The coefficients $\L _n$ computed by (\ref{EKL})--(\ref{BE}) up to $n=4000$, 
on a logarithmic $n$-scale (straight line: the function $(\log n +C)$ of (\ref{LRH})).}
\end{figure}

Low-$n$ calculations of $\L _n$ (fig. 3) agree very early with the \emph{logarithmic} behavior (\ref{LRH}),
just as they agreed for $\l _n$ with its leading behavior under RH \cite{K}\cite{M1}. 
The remainder term $\delta \L _n  = \L _n -(\log n + C)$ looks compatible with an $o(1)$ bound (fig.~4), 
albeit much less neatly than $\delta \l ^{\rm K}_n$ \cite[fig.~1]{K}\cite[fig.~6b]{M1},
(note: both of these plot $\delta \l ^{\rm L}_n = n \ \delta \l ^{\rm K}_n $). For the record,
\beq
\textstyle \L _1 = \frac{3}{2} \log \frac{\pi}{3} \approx 0.069176395771 ,
\ \ \L _2 \approx 0.22745427267 , \ \ \L _3 \approx 0.45671413349;
\eeq
\vskip -39pt

\bea
\L _{10000} \si\approx\se 8.428662659671506 \qquad (\delta \L _{10000} \approx +0.0020794), \nonumber \\[-2pt]
\L _{20000} \si\approx\se 9.120189975922122 \qquad (\delta \L _{20000} \approx -0.000485565), \nonumber
\eea
It would be interesting to comprehend the bumpy fine structure of $\delta \L _n$.

\begin{figure}
\center
\includegraphics[scale=.4]{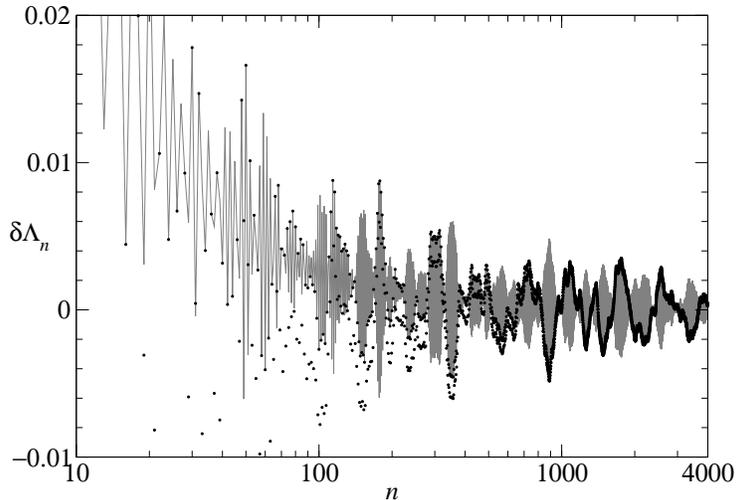}
\caption{\small The remainder sequence $\delta \L _n = \L _n -(\log n+C)$ 
(in gray: the connecting segments are drawn for visual clarity only), 
and a rectified form $(-1)^n \delta \L _n$ (black dots) to cancel the period-2 oscillations.}
\end{figure}

\subsection{Imprints of putative zeros violating RH}

RH-violating zeros $\rho $ (if any) seem to enter the picture just as for the $\l _n$:
their contributions (\ref{ST1}) will asymptotically dominate $\log n$,
but numerically they will emerge and take over extremely late. 
For such a zero $\rho = \hf +t+ \mi T$, with $0<t<\hf$ and $T \gtrsim 2.4 \cdot 10^{12}$ \cite{G},
its contribution sizes like $T^{-2}(2n/T)^t /\log n$ in modulus, by (\ref{LABS}).
We then get its crossover threshold (in order of magnitude, neglecting logarithms and constants relative to powers) 
by solving 
\bea
\label{TN}
T^{-2}(n/T)^t \si\approx\se 1 \\
\Longrightarrow \qquad n \si\gtrsim\se T ^{1+2/t} \qquad 
(\mbox{best case: }O(T^{5+\eps }) \mbox{ for } t=\hf-0) .
\eea
This is worse than (\ref{LIt}) for $\l _n$, all the more if a negativity test is pursued
(the right-hand side of (\ref{TN}) must then be $\log ^2 n$). 
There is however room for possible improvement:
the core problem is to filter out a weak $\rho $-signal from the given background (\ref{LRH}),
therefore any predictable structure in the latter is liable to boost the gain.
For instance, the hyperfine structure of $\delta \L _n$ is oscillatory of period 2 (fig.~4);
this suggests to average over that period, which \emph{empirically}
discloses a rather neat $(1/n)$-decay trend (fig.~5):

\begin{figure}
\center
\includegraphics[scale=.4]{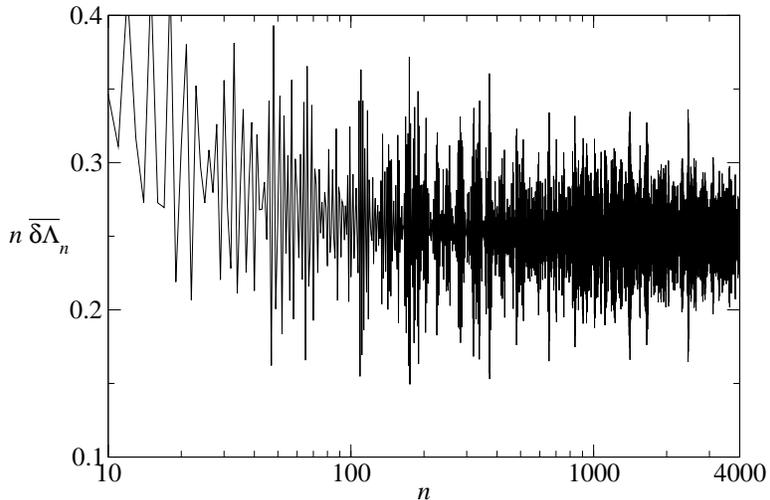}
\caption{\small The averaged remainder sequence (\ref{LAV}) rescaled by $n$, namely: $n \, \overline{\delta \L }_n$.
(Some further values: $0.27027$ for $n=10000,\ 0.23970$ for $n=20000$.)
}
\end{figure}

\beq
\label{LAV}
\overline{\delta \L }_n \defi \hf (\delta \L _n + \delta \L _{n-1}) \approx 0.25 /n .
\eeq
The same operation on a $\rho $-signal $F_n (\rho )$ in (\ref{LNRH}) roughly applies $\hf (\d /\d n)$
to the factor $n^T$ therein (again neglecting $t \ll T$ and $\log n$), i.e., multiplies it by $\hf (T/n)$.
Thus heuristically, i.e., conjecturing the truth of (\ref{LAV}) for $n \to \infty$ under RH, 
the crossover condition improves from (\ref{TN}) to
\bea
\label{TNI}
(T/n) \, T^{-2}(n/T)^t \si\approx\se \overline{\delta \L }_n \approx 1 /n \nonumber \\
\Longrightarrow \quad n \si \gtrsim\se T^{1+1/t} 
\ (\mbox{best case: }O(T^{3+\eps }) \mbox{ for } t=\hf-0) . \quad
\eea
We can hope that efficient signal-analysis techniques may still lower this detection threshold.
And an empirical attitude may suffice here: once a violating zero would be suspected and roughly located,
other rigorous algorithms exist to find it accurately (or disprove it). \cite{G}

\subsection{\emph{The} hitch}

A major computational issue is that, according to (\ref{EKL}), the $(\log n)$-sized values $\L _n$ 
result from alternating summations of \emph{much faster-growing} terms: 
this entails a loss of precision increasing with $n$. 
Thus in our case (sums $\sum s_m$ of order comparable to unity), to reach the \emph{slightest} end accuracy
we must use each summand $s_m$ up to $\approx \log_{10} |s_m|$ significant digits (in base 10 throughout);
plus uniformly $D$ \emph{more} to obtain $\sum s_m$ accurate to $D$ digits.

We quantify the precision loss in (\ref{EKL}) at large fixed $n$ by using the Stirling formula,
to find that $m_\ast \approx n/\sqrt 2$ is where the largest summand occurs
and the minimum required precision $\log_{10} |s_m|$ peaks, reaching
$\log_{10} |A_{nm_\ast } \log 2 \xi (2m_\ast )|$ $\sim \log_{10} (3+2 \sqrt 2) \, n \approx 0.76555 \, n$ 
digits, see fig.~6 (vs a precision loss $\approx 0.25 \, n$ digits for $\l _n$ \cite[fig.~6]{M2}).
Even then, a crude feed of (\ref{EKL}), (\ref{Udef}) or (\ref{Vdef}) 
into a mainstream arbitrary-precision system (Mathematica 10 \cite{W}) 
suffices to readily output the $\L _n$-values of \S~4.1.
Computing times varied erratically but could go down to ca. 4 min for $\L _{10000}$, 43 min for $\L _{20000}$ 
using (\ref{Vdef}) (CPU times on an Intel Xeon E5-2670 0 @ 2.6 GHz processor).
\vskip 6.5mm

\begin{figure}[h]
\center
\includegraphics[scale=.4]{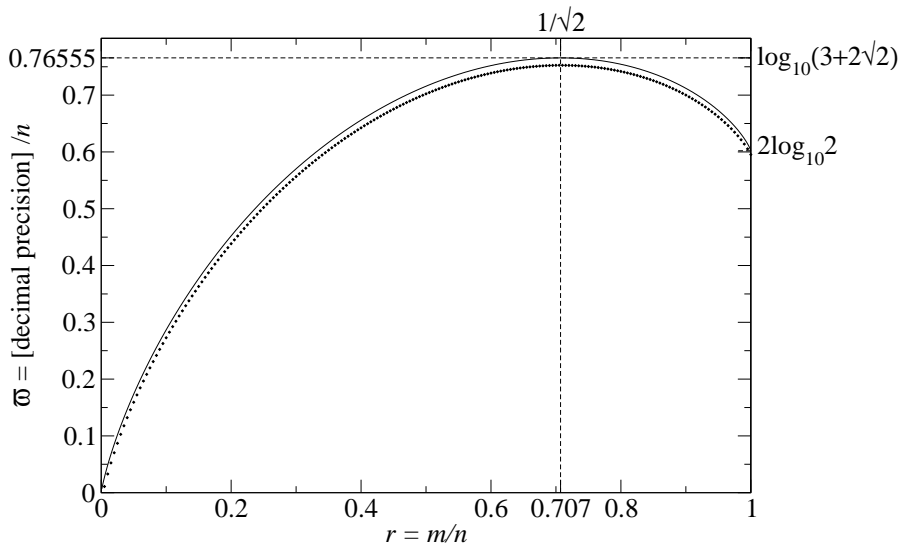}
\caption{\small Minimum decimal precisions needed for the summands of $\L _n$ in (\ref{EKL}),
as estimated by $\log_{10} |A_{nm} \log 2\xi (2m)|$ which is plotted against $m$ in axes rescaled by $1/n$.
Dotted curve: the case $n=200$; continuous curve: the $n \to \infty$ limiting form 
$\varpi = -2 \, r \log_{10} r + (1+r) \log_{10} (1+r) -(1-r) \log_{10} (1-r) \quad (r=m/n)$.
}
\end{figure}

Now with $|T|\gtrsim 2.4 \cdot 10^{12}$ currently, the challenge is to probe $n \gtrsim 2 \cdot 10^{36}$
(if the more favorable estimate (\ref{TNI}) holds, $10^{60}$ otherwise), 
which then needs a working precision $\gtrsim 1.6 \cdot 10^{36}$ decimal places at times.
This need for a huge precision already burdened the original $\l _n$ 
but somewhat less and amidst several steeper complexities, 
now for the $\L _n$ the ill-conditioning worsened while the other difficulties waned.

As advantages of $\{\L _n\}$ over $\{\l _n\}$, inversely: the $\L _n$ are fully explicit;
their evaluations are recursion-free, thus \emph{very few} samples (at high enough~$n$, for sure)
might suffice to signal that RH is violated \emph{somewhere};
and the required working precision peaking at $\approx 0.766 \, n$ stands as the \emph{only} stumbling block,
and as a \emph{purely logistic} problem, which might still be eased if (\ref{EKL}) came to admit better conditioned variants.
Thus in (\ref{Udef}), a much lower precision (growing like $\hf \log_{10} n$)
suffices for $\log 2\pi $ with its factor $(2A_{n0})^{-1} \sim -\sqrt{\pi n}/2$ which grows negligibly, 
compared to the ${A_{nm} \log(|B_{2m}|/(2m-3)!!) :}$ only these \emph{simpler} expressions demand 
maximal precision, and only for $m \approx n/ \sqrt 2$.
\bigskip

While other sequences sensitive to RH for large $n$ are known \cite{BD}\cite{FV}, not to mention Keiper--Li again,
we are unaware of any previous case combining a fully \emph{closed form} like (\ref{EKL})
with a practical sensitivity-threshold of \emph{tempered growth} $n=O(T^\nu )$.

\begin{appendix}
\section*{Appendix: Centered variant}

We sketch a treatment parallel to the main text for our Li-type sequences 
having the alternative basepoint $x=\hf$ (the \emph{center} for the $\xi $-function).

We recall that the Functional Equation $\xi (1-x) \equiv \xi (x)$ allows us,
in place of the mapping $z \mapsto x=(1-z)^{-1}$ within $\xi $ as in (\ref{KEd}), to use 
the double-valued one $y \mapsto x_{\tilde w} (y)= \hf \pm \sqrt{\tilde w} \, y^{1/2} / (1-y)$ 
on the unit disk (parametrized by $\tilde w>0$).
That still maps the unit circle ${\{ |y|=1 \}}$ to the completed critical line $L \cup \infty$, 
but now minus its interval $\{ |\Im x|< \hf \sqrt{\tilde w} \}$.
As before, all Riemann zeros on $L$ have to pull back to ${\{ |y|=1 \}}$
which then imposes $\tilde w < 4 \min_\rho |\Im \rho |^2 \approx 799.1618$. 
We thus define the sequence $\{ \l ^0_n (\tilde w) \}$ by
\beq
\label{Cdef}
\log 2\xi \Biggl( \hf \pm \frac{\sqrt{\tilde w} \, y^{1/2}}{1-y} \Biggr) \equiv 
\log 2\xi (\hf ) + \sum_{n=1}^\infty \frac{\l ^0_n (\tilde w)}{n} \, y^n 
\eeq
(\cite[\S 3.4]{VK}, where only the case $\tilde w=1 $ is detailed), or
\beq
\label{CRs}
\frac{\l ^0_n (\tilde w)}{n} \equiv \frac{1}{2 \pi\mi} \oint \frac{\d y}{y^{n+1}} 
\log 2\xi \bigl( x_{\tilde w} (y) \bigr) , \qquad n=1,2,\ldots 
\eeq

We now build an \emph{explicit} variant for this sequence (\ref{CRs}),
similar to $\{ \L _n \}$ for $\{ \l _n^{\rm K} \}$.
First, the deformations of (\ref{CRs}) analogous to those in \S 2.1 read as
\beq
\frac{1}{2 \pi\mi} \oint \frac{\d y}{B_{y_0}(y) \cdots B_{y_n}(y)} 
\log 2\xi (x)  \qquad (\mbox{here } x \equiv  x_{\tilde w} (y)) ,
\eeq
for which the simplest analytical form we found, similar to (\ref{KLD}), is now 
\beq
\label{CR}
\frac{1}{2 \pi\mi} \oint \frac{2 \, \d r}{(r+1)^2} \prod_{m=0}^n \frac{r+r_m}{r-r_m} 
\, \log 2\xi (x) , \qquad r_m \defi \frac{1+y_m}{1-y_m} ,
\eeq
all in terms of the new variable 
\beq
r \defi \frac{1+y}{1-y} \equiv \bigl[ 1+(2x-1)^2/{\tilde w} \bigr] ^{1/2} \qquad (\Re r>0) .
\eeq
Then with $x_m \equiv 2m$ as before (but now including $m=0$),
the integral (\ref{CR}) evaluated by the residue theorem yields the explicit result 
(akin to (\ref{EKL})--(\ref{BE}))
\beq
\L ^0_n (\tilde w) \defi \sum_{m=1}^n \frac{2}{(r_m \!+\! 1)^2}
\frac{\prod\limits_{k=0}^n (r_m \!+\! r_k)}{\prod\limits_{k \ne m} (r_m\!-\! r_k)} \log 2\xi (2m) , 
\quad r_m \equiv \sqrt{ 1+(4m \!-\! 1)^2/{\tilde w}} .
\eeq
This result is, however, algebraically less simple and less analyzable than for $\L _n$ before.
A potential asset is that it openly relies on the Functional Equation,
but we saw no practical benefit accruing from that yet.

The corresponding \emph{asymptotic alternative} for RH analogous to (\ref{LNRH})--(\ref{LRH}) reads as
\bea
\label{CNRH}
\bullet \ \mbox{RH false:} \ \L _n \si\sim\se 
\Biggl\{ \sum_{\Re \rho >1/2} \Delta_\rho \L_n^0 (\tilde w) \Biggr\}
\pmod {o(n^\eps)\ \forall \eps >0} \\
\label{ST0}
\si\se \mbox{with } \log |\Delta_\rho \L_n^0 (\tilde w)| \sim (\rho  -1/2) \log n , \\
\label{CRH}
\bullet \ \mbox{RH true:} \ \L _n \si\sim\se \sqrt{\tilde w} \, (\log n + C) , \quad 
C = \hf(\g \!-\! \log \pi \!-\! 1) \approx -0.78375711 . \qquad
\eea
The latter is proved by extending Oesterl\'e's method just as with $\L _n$; whereas
the former needs large-$n$ estimations of the product in (\ref{CR}), 
but the ones we have remain crude compared to the full Stirling formula available for (\ref{GN});
that precludes us from reaching the absolute scales of the $\Delta_\rho \L_n^0 (\tilde w)$
and the values of $n$ from which any such terms might become detectable.

As for numerical tests, all results are very close to those shown above for $\L _n$,
aside from the overall factor $\sqrt{\tilde w}$ in (\ref{CRH}) 
(but nothing about the case RH false can be tested: that is still way beyond numerical reach).
\end{appendix}

\end{document}